\documentclass[12p]{article}
\usepackage[dvips]{graphicx}
\usepackage{subfig}
\usepackage{multirow}

\begin{document}

\title {A sharp bound on the expected local time of a continuous ${\cal L}_2$-bounded Martingale}

\author{
David Gilat and Isaac Meilijson
\\
School of Mathematical Sciences, R. and B. Sackler Faculty of Exact Sciences \\
Tel Aviv University, Tel Aviv 6997801, Israel
\\
{\em E-mail: \tt{isaco@tauex.tau.ac.il, gilat@tauex.tau.ac.il}}
\\
Laura Sacerdote
\\
Department of Mathematics G. Peano,	University of Torino \\
	Via Carlo Alberto 10, 10123 Torino, Italia \\
{\em E-mail: \tt{laura.sacerdote@unito.it}}
}

\maketitle

\pagenumbering{arabic}

\begin{abstract}
For a continuous ${\cal L}_2$-bounded Martingale with no intervals of constancy, starting at $0$ and having final variance $\sigma^2$, the expected local time at $x \in \cal{R}$ is at most $\sqrt{\sigma^2+x^2}-|x|$. This sharp bound is attained by Standard Brownian Motion stopped at the first exit time from the interval $(x-\sqrt{\sigma^2+x^2},x+\sqrt{\sigma^2+x^2})$. Sharp bounds for the expected maximum, maximal absolute value, maximal diameter and maximal number of upcrossings of intervals, have been established by Dubins and Schwarz (1988), Dubins, Gilat and Meilijson (2009) and by the authors (2017).


 
\end{abstract}

\section{Introduction}

For ${\cal L}_2$-bounded Martingales with mean zero and final variance $\sigma^2$, sharp upper bounds for the expectations of the maximum $M$ (the bound is $\sigma$), maximal absolute value $|M|$ (the bound is $\sigma \sqrt{2}$), diameter or range $D=M-m$ (where $m$ is the minimum, and the bound is $\sigma \sqrt{3}$) and number of up-crossings of a given interval $(x,b)$, were established in Dubins \& Schwarz \cite{DS2} for $M$ and $|M|$, in Dubins, Gilat \& Meilijson \cite{DGM} for $D$ and in Gilat, Meilijson \& Sacerdote \cite{GMS} for up-crossings. It seems only natural for {\em local time} at a given $x$ to be next in line. While all the former quantities are well defined for all (c\'{a}dl\`{a}g version) Martingales, the extent of the definition of local time and its meaning as density of the occupation measure, have not been fully studied. Attention is therefore focused on stopped Brownian Motion (Sections 2 and 3) and subsequently (Section 5) restricted to continuous Martingales to which the Brownian Motion result easily extends by embedding \cite{Bjork, Doeblin, DS1, Monroe}. The bound for local time at $x$ in this class of martingales is $\sigma [\sqrt{1+({x \over \sigma})^2}-{|x| \over \sigma}]$. This is a decreasing function of $|x|$, with maximal value $\sigma$ at $x=0$, asymptotic to ${{\sigma^2} \over {2 |x|}}$ for large $|x|$. For the case where $b-x \le \sigma$, the sharp bound for up-crossings (\cite{GMS}) is ${\sigma \over {2 (b-x)}}[\sqrt{1+({x \over \sigma})^2}-{|x| \over \sigma}]$. Thus, the sharp bound for local time is consistent with a definition of local time in terms of up-crossings \cite{MortPers}.

\section{The case of stopped Brownian Motion}

\noindent {\bf A remark on the role of boundedness in ${\cal L}_2$.} Let $\{B(t) ; t \ge 0 \ , \ B(0)=0\}$ be standard Brownian Motion and let $\tau$ be a stopping time on $B$ such that the Martingale $W(t)=B(\tau \wedge t)$ is bounded in ${\cal L}_2$. Let $\tau_t = \tau \wedge t$. Since $B(t)^2-t$ is a mean-zero Martingale, $E[\tau_t]=E[W(t)^2]$ is uniformly bounded. Hence, by Lebesgue's monotone convergence theorem, $E[\tau]<\infty$. But then the Martingale property of $B(t)^2 - t$ extends to stopping at $\tau$, yielding $E[\tau]=E[B(\tau)^2]=\sigma^2$.

It is worth noticing that finite final variance is weaker than boundedness in ${\cal L}_2$. Let $F$ be a distribution with mean zero and finite variance, and let $\tau$ be the first hitting time of $B$ to $X$, where $X$ had been independently drawn from $F$. Then $B(\tau)$ is distributed $F$ and $W$ has the "correct" final variance, but it is not a Martingale bounded in ${\cal L}_2$.

\bigskip

With this clarification in mind, let $\{B(t) ; t \ge 0 \ , \ B(0)=0\}$ be SBM, let $\sigma \in (0,\infty)$ and let $\tau$ be a stopping time on $B$ with $E[B(\tau)^2]=E[\tau]=\sigma^2$. Let $x \in \cal{R}$. This section introduces a stopping time that maximizes the expected local time $E[L_x(\tau)]$ at $x$.

Local time for Brownian Motion was first introduced by L\'{e}vy \cite{Levy}. It has been equivalently defined as the density of an occupation measure or directly as a limit of occupation time or up-crossings of a vanishing interval \cite{KarShr, RogWill}, as well as via its role in the explicit submartingale Doob-Meyer \cite{Doob, Meyer1, Meyer2} decomposition of absolute Brownian Motion as provided by Tanaka's formula \cite{KarShr}.

For $x=0$ this decomposition becomes $|B(\tau)|=L_0(\tau) + S(\tau)$, where $S$ is a mean-zero Martingale, so $E[L_0(\tau)]=E[|B(\tau)|]$ is the ${ \cal L}_1$ norm of a random variable $B(\tau)$ with mean $0$ and ${\cal L}_2$ norm $\sigma$. Since the former is bounded from above by the latter and can be made equal to it by letting $\tau$ be the first exit time by $B$ from the interval $(-\sigma, \sigma)$, this is clearly an optimal solution, and the supremum of $E[L_0(\tau)]$ over all feasible $\tau$ is $\sigma$.

An alternative proof can be based on Dubins \& Schwarz \cite{DS2} sharp upper bound $\sigma$ on the expected maximum of stopped SBM, via Levy \cite{Levy, MortPers} celebrated joint distribution equality of the two pairs of processes $(|B|,L_0)$ and $(M-B,M)$, where $M$ is the cumulative maximum process of $B$, $M(t)=\max \{B(s) \ | \ s \le t\}$.

For the case $x \neq 0$, assume without loss of generality that $x>0$, and consider any stopping time $\tau$ with $E[\tau]=\sigma^2$. Let $T_x$ be the time of the first visit by $B$ to $x$ and let $A=\{T_x < \tau\}$ be the event that $B$ visits $x$ before time $\tau$. $P(A)<1$ since $E[T_x]=\infty$. As clearly $\sup_{\tau}E[L_x(\tau)]>0$, we may assume that $P(A)>0$, because otherwise $L_x(\tau)=0$ a.s. Since
\begin{eqnarray} \label{EBtau}
0 & = & E[B(\tau)] \label{EBtau} \\
& = & P(A) E[B(\tau) | A]+(1-P(A)) E[B(\tau) | A^c]=P(A) x + (1-P(A)) E[B(\tau) | A^c] \nonumber
\end{eqnarray}
denote $y=E[B(\tau) | A^c]=-x {{P(A)} \over {1-P(A)}}$,
let $v=\mbox{Var}[B(\tau) | A^c]$ and $\eta^2 = \mbox{Var}[B(\tau) | A] = {{\sigma^2 + x y - (1-P(A)) v} \over {P(A)}}$.

The choice of $\tau$ on the event $A$ is constrained only by the remaining variance budget $\eta^2$. By the above analysis for $x=0$, the expected local time at $x$ can be increased if $\tau$ is replaced on $A$ by the first exit time after time $T_x$ from the interval $(x-\eta, x+\eta)$. This yields conditional expected local time $\eta$ and unconditional expected local time $\eta P(A)$. Hence, for fixed $P(A)$, the expected local time at $x$ can be further increased by raising $\eta$ as much as possible. This is achieved by rendering $B(\tau)$ a.s. constant off $A$, i.e., by forcing $v=0$. In other words, the search for $\tau$ can be restricted to those of the form: Initially wait until the first exit time by $B$ from the interval $(y,x)$. If at $y$, full stop. If at $x$ (and this happens with probability $P(A)={{-y} \over {x-y}}$), stop at the first exit time by $B$ after time $T_x$ from the interval $(x-\eta, x+\eta)$, where $\eta^2 = \mbox{Var}[B(\tau) | A] = {{\sigma^2 + x y} \over {P(A)}} = {{(\sigma^2 + x y)(x-y)} \over {-y}}$.

It only remains to maximize over $y<0$ the objective function $\eta P(A)=\sqrt{{{-(\sigma^2+ x y) y} \ \over {x-y}}}$. As can be easily calculated, this is achieved by $y=x-\sqrt{\sigma^2+x^2}$. Under this choice of $y$, $\eta=\sqrt{\sigma^2+x^2}$ and the maximized objective value is $\sqrt{{{-(\sigma^2+ x y) y} \ \over {x-y}}} = \sqrt{\sigma^2+x^2}-x$. As a consequence, since $x-\eta=y$, the ensuing $\tau$ is simply the first exit time by $B$ from the interval $(x-\sqrt{\sigma^2+x^2}, x+\sqrt{\sigma^2+x^2})$.

The bound $\sqrt{\sigma^2+x^2}-x$, that coincides with $\sigma$ for $x=0$, is conveniently expressed as
\begin{equation} \label{oneoverx}
(\sqrt{\sigma^2+x^2}-x){{\sqrt{\sigma^2+x^2}+x} \over {\sqrt{\sigma^2+x^2}+x}} = {{\sigma^2} \over {\sqrt{\sigma^2+x^2}+x}} \approx {{\sigma^2} \over {2 x}}
\end{equation}
to see its behavior for large $x$.

\section{An invariance property of expected local time of stopped Brownian Motion}

It is shown in this section that the expected local time $E[L_x(\tau)]$ at $x$ is fully determined by the distribution of the r.v. $X=B(\tau)$. In fact,
\begin{eqnarray} 
E[L_x(\tau)] & = & E[|X-x|]-|x|  \nonumber \\
& = & 2 E[(X-x)^+] = 2 \int_x^\infty P(X>t)dt \ ; \ x \ge 0 \nonumber \\
& = & 2 E[(X-x)^-] = 2 \int_{-^\infty}^x P(X<t)dt  \ ; \ x \le 0 \label{ELX}
\end{eqnarray}

To prove (\ref{ELX}), apply on the event $A$ the Doob-Meyer decomposition of the previous section to the mean-zero Brownian motion $W(t)=B(T_x + t)-x$:
\begin{equation} \label{IAB}
I_A |B(\tau)-x| = I_A (L_x(\tau) + S(\tau)) = I_A L_x(\tau) + I_A S(\tau) = L_x(\tau) + I_A S(\tau)
\end{equation}
where $S$ is a mean-zero Martingale with respect to a filtration whose starting $\sigma$-field contains $A$. As this entails $E[I_A S(\tau)]=0$, it follows for $x \ge 0$ (similarly for $x \le 0$) that
\begin{equation} \label{ELXtau}
E[L_x(\tau)]=E[I_A |B(\tau)-x|]=2 E[I_A (B(\tau)-x)^+]=2 E[(B(\tau)-x)^+] \equiv 2 E[(X-x)^+]
\end{equation}
where the third equality is due to the fact that $\{B(\tau)>x\} \subseteq A$.

To illustrate this formula with the dichotomous result of the previous section, $2 E[(X-x)^+] = 2 \sqrt{\sigma^2+x^2}
{{\sqrt{\sigma^2+x^2}-x} \over {2 \sqrt{\sigma^2 + x^2}}} = \sqrt{\sigma^2+x^2}-x$. It is easy to check that this dichotomous distribution maximizes $E[|X-x|]$ (i.e., $|x|+E[L_x(\tau)]$) among all distributions (of $X=B(\tau)$) with mean zero and variance $\sigma^2$, thus providing yet another proof of our basic result.

For $X=B(\tau)$, $\tau$ the first exit time from an interval $I=(-a,b)$ containing $0$, $SD[X]=\sigma=\sqrt{a b}$. The expected local time at $x \in I$ until time $\tau$ is ${{2 a (b-x)} \over {a+b}}=\sigma {2 \over {\sqrt{{a \over b}}+\sqrt{{b \over a}}}}(1-{x \over b})$ for $x \in [0,b]$ and ${{2 b (a+x)} \over {a+b}}=\sigma {2 \over {\sqrt{{a \over b}}+\sqrt{{b \over a}}}}(1-{{|x|} \over a})$ for $x \in [-a,0]$. For $x=0$, this is the harmonic mean of $a$ and $b$.

If $X$ is normally distributed, expected local time at $x$, $2 E[(X-x)^+]= 2 E[(\sigma Z-x)^+]= 2 \int_x^\infty (\sigma z - x) \phi(z)dz= 2 \sigma \phi(x) - 2 x \Phi^*(x)$, is $\sigma\sqrt{{2 \over \pi}}$ at $x=0$ (compared to the upper bound $\sigma$), and decreases to zero faster than the density $\phi$, as compared to the rate ${1 \over x}$ of the bound.

If $X=\sigma (Y-1)$ for $Y \sim Exp(1)$, expected local time $\sigma {2 \over e} \exp\{-{x \over \sigma}\}$ at $x \ge 0$ is indeed below $\sigma$ and decays exponentially.

\noindent {\bf Remark}: It is clear from the RHS of (\ref{ELX}) that the expected local time of SBM at $x$ until any integrable stopping time is unimodal in $x$, with a maximum at $x=0$. It appears that the certainty of visits at the initial point takes precedence under any terminal distribution.

\section{A digression on the Chacon-Walsh solutions to Skorokhod's embedding problem: Further roles of the function $E[|X-x|]$}

Skorokhod \cite{Skorokhod} posed the question of existence of (in our context, integrable) stopping times $\tau$ on SBM $B$ with arbitrary mean-zero, finite variance distributions for $B(\tau)$. He solved the problem with randomized stopping times, first exit times from intervals, by representing the distribution as a mixture of dichotomous mean-zero distributions. Of the many non-randomized solutions, we focus on the Chacon-Walsh \cite{ChaWal} idea, based on representing the r.v. $X$ whose distribution is to be embedded in SBM, as the limit of a Martingale with dichotomous transitions. The embedding stopping time is thus the limit of a sequence of first-exit times from intervals. The Chacon-Walsh method represents the function $E[|X-x|]$ as the envelope of a sequence of linear supporting minorants, $|x|$ being the function corresponding to the point mass at zero, that plays the role of initial distribution. It is a celebrated theorem of Hardy, Littlewood \& Polya (\cite{HardyLitt}, extended and re-analyzed by Strassen \cite{Strassen}, Rothschild \& Stiglitz \cite{RothStig} and others) that a distribution $G$ is a Martingale dilation of another, $F$ ($(F,G)$ are the respective marginals of some Martingale pair $(X,Y)$) if and only if $E[|Y-x|] \ge E[|X-x|]$ pointwise. The Chacon-Walsh technique has been extended to arbitrary initial distributions, i.e., to embed a pair of distributions related by the inequality above, with $B(\tau_1) \sim F$, $B(\tau_2) \sim G$ and $\tau_1 \le \tau_2$ a.s. This construction provides a didactic proof of the Hardy-Littlewood-Polya characterization of Martingale dilation. The common theme of this notion and the subject matter of this study and its direct predecessors is the general question of how does the final distribution of a Martingale restrict its overall variability.

\section{Expected local time of more general Martingales}

C\'{a}dl\`{a}g Martingales are known to be optional sampling of Brownian Motion (Monroe \cite{Monroe}), i.e., can be represented as $SBM$ sampled at an increasing family $\{\tau_t \ ; \ t \ge 0\}$ of stopping times. If the Martingale has a.s. continuous paths and no intervals of constancy, this time-change can be taken to be continuous (Doeblin \cite{Doeblin}, Dubins and Schwarz \cite{DS1}). Bj\"{o}rk \cite{Bjork} (Proposition 4.1) has shown that such continuous adapted time changes preserve local time. Hence, the bounds on expected local time in Brownian Motion extend verbatim to such Martingales.

It may be possible to extend the scope of the bounds presented above to more general Martingales, but this would require a conceptually meaningful definition of local time, hopefully preserving its nature as occupation measure density. It makes sense to expect that discontinuities reduce local time relative to that of the underlying Brownian Motion (so upper bounds would hold a fortiori) but intervals of constancy present a challenge:

Consider a continuous time Markov chain starting at zero that changes its value by equally likely $\pm 1$  at the consecutive arrival points of a Poisson process. This Markovian Martingale has sojourns with positive Lebesgue measure at the integers, so the regular definition of local time makes local time infinite as soon as it is positive. This defies any bound as claimed above. Perhaps some variant of the definition of local time, such as the alternative definition in terms of upcrossings (e.g., M\"{o}rters and Peres \cite{MortPers}) could accommodate these pathologies.

\section*{Acknowledgements}

We thank Ioannis Karatzas for useful discussions. This work was motivated by his lecture in a memorial conference for Larry Shepp organized by Philip Ernst at Rice University, 2018. Laura Sacerdote's research is partially supported by INDAM - GNCS and by a University of Torino 2019 grant for Stochastic and Statistical Methods and Models. Isaac Meilijson research is partially supported by grant 217/16 of the Israel Science Foundation.

\end{document}